\documentclass[11pt]{article}
\usepackage{latexsym, amsmath, amssymb}

\usepackage{amsmath}
\usepackage{amsfonts}
\usepackage{amssymb}
\usepackage{amscd}
\usepackage{amsthm}
\usepackage{fancyhdr}
\theoremstyle{plain}
\newtheorem{theorem}{Theorem}[section]
\newtheorem{proposition}[theorem]{Proposition}
\newtheorem{remark}[theorem]{Remark}
\newtheorem{lemma}[theorem]{Lemma}
\newtheorem{definition}[theorem]{Definition}
\newtheorem{corollary}[theorem]{Corollary}

\newcommand\sO{{\mathcal O}}

  \def \tab#1{\kern #1 truein}

\begin{document}
 \title{Regularity and Cohomological Splitting Conditions for Vector Bundles on Multiprojective Spaces}
\author{Edoardo Ballico and Francesco Malaspina
\vspace{4pt}\\
{\small   Universit\`a di Trento}\\
{\small\it 38123 Povo (TN), Italy}\\
{\small\it e-mail: ballico@science.unitn.it}\\
\vspace{6pt}\\
         {\small  Politecnico di Torino}\\
{\small\it  Corso Duca degli Abruzzi 24, 10129 Torino, Italy}\\
{\small\it e-mail: francesco.malaspina@polito.it}} \maketitle \def\thefootnote{}
\footnote{
keywords: vector bundles; multiprojective spaces; Castelnuovo-Mumford regularity.}
  \begin{abstract}We give a definition of regularity on multiprojective spaces  which is
different from the definitions of Hoffmann-Wang and  Costa-Mir\'o-Roig. Using this notion we prove some splitting criteria
for vector bundles.
  \end{abstract}

   \section{Introduction}\label{S1}
    In chapter $14$ of \cite{m} Mumford introduced the concept of regularity for a coherent sheaf on a projective space $\bf P^{n}$.
It was soon clear that Mumford's definition of Castelnuovo-Mumford regularity
was a key notion and a fundamental tool in many areas of algebraic geometry and commutative algebra. Intuitively, it measures the complexity of a module or sheaf. The regularity of a coherent sheaf estimates the smallest twist for which the sheaf is generated by its global sections. Several extensions
of this notion have been proposed over the years to handle different situations (\cite{am}, \cite{bm2}, \cite{c}, \cite{cm2}, \cite{CM3}, \cite{ms}, \cite{sv} and \cite{hw}).
Maclagan and  Smith   developed in \cite{ms} an interesting variant of multigraded Castelnuovo-Mumford regularity, motivated by toric geometry.
 For a different approach to multigraded Castelnuovo-Mumford regularity, see \cite{sv}.
 Hoffmann and Wang in \cite{hw} gave the following definition of regularity
on $\mathbb {P}^n\times \mathbb {P}^m$:
\begin{definition}[Hoffmann and Wang]\label{d0}
A  coherent sheaf $F$ on $X$  is said to be \emph{HW-$(p,p')$-regular} or $(p,p')$-regular in the sense of
Hoffmann and Wang if, for all $i>0$, $$H^i(F(p,p')\otimes
\sO(j,k))=0$$
whenever $j+k=-i-1$, $j< 0$ and $ k< 0$.
We will say ``~\emph{HW-regular}~'' instead of ``~HW-$(0,0)$-regular~'',
and ``~\emph{HW-$p$-regular}~'' instead of ``~HW-$(p,p)$-regular~''.
We define the {\it HW-regularity} of $F$, $\mbox{HW-Reg}(F)$, as the smallest integer $p$ such that $F$ is HW-$p$-regular. We set
$\mbox{HW-Reg}(F)=-\infty$ if there is no such integer.
\end{definition}
For a definition of regularity on multiprojective spaces using Beilinson's type spectral sequence,
see \cite{cm2}. Here we introduce the following modification
of Hoffman and Wang's definition:
\begin{definition}\label{d1}
A  coherent sheaf $F$ on $X$  is said to be  \emph{$(p,p')$-regular} if, for all $i>0$, $$H^i(F(p,p')\otimes
\sO(j,k))=0$$
whenever $j+k=-i$, $-n\leq j\leq 0$ and $-m\leq k\leq 0$.\\
We often say ``~\emph{regular}~'' instead of ``~$(0,0)$-regular~'',
and ``~\emph{$p$-regular}~'' instead of ``~$(p,p)$-regular~''.
We define the {\it regularity} of $F$, $Reg (F)$, as the least integer $p$ such that $F$ is $p$-regular. We set
$Reg (F)=-\infty$ if there is no such integer.
\end{definition}
Our motivation for studying regularity over multigraded polinomial rings is  to prove splitting criteria for vector bundles on multiprojective spaces. A well known result by Horrocks (see \cite{Ho}) characterizes  vector bundles without intermediate cohomology on a projective space as  direct sums of line bundles.
This criterion fails on  more general varieties. In fact there exist non-split vector bundles  without intermediate cohomology, called ACM bundles. The theory of regularity  allows to prove easily Horrocks criterion on $\mathbb {P}^n$ and its improvement by
Evans-Griffith depending on the rank of the vector bundle (see \cite{e}).
In \cite{bm2} we introduced the notion of Qregularity on a quadric hypersurface in order to prove an extension of Evans-Griffith criterion to vector bundles on smooth quadric hypersurfaces. In particular we obtained a new and simple proof of the Kn\"{o}rrer's characterization of ACM bundles.

 In section \ref{S2} we will prove that our definition of regularity for
biprojective spaces  satisfies the analogues of the classical properties on $\mathbb P^{n}$. Moreover it has several nice features and allows us to classify some ``~extremal cases~''.

 In section \ref{S3} we will apply our definition of regularity to prove a few splitting criteria for vector bundles on $\mathbb {P}^n\times \mathbb {P}^m$.
On $\mathbb {P}^n$ all  line bundles are ACM but on $\mathbb {P}^n\times \mathbb {P}^m$ there are line bundles which are not ACM.
We prove the following extension of the Horrocks criterion on $\mathbb {P}^n\times \mathbb {P}^m$:
\begin{theorem}\label{t1}Let $E$ be a rank $r$ vector bundle on $\mathbb {P}^n\times \mathbb {P}^m$. Then the following conditions are equivalent:
  \begin{enumerate}
  \item[(i)] for any $i=1, \dots, m+n-1$ and for any integer $t$,  $H^i(E(t,t)\otimes \sO(j,k))=0$ whenever
  $j+k=-i$, $-n\leq j\leq 0$ and $-m\leq k\leq 0$.
  \item[(ii)] There are $r$ integers $t_1, \dots, t_r$ such that $E\cong \bigoplus_{i=1}^r \sO(t_i,t_i)$.
  \end{enumerate}
  \end{theorem}

Then we prove the following result:
  \begin{theorem}\label{t2}Let $E$ be a  vector bundle on $\mathbb {P}^n\times \mathbb {P}^m$. Then the following conditions are equivalent:
  \begin{enumerate}
  \item[(a)] for any $i=1, \dots, m+n-1$ and for any integer $t$,  $H^i(E(t,t)\otimes \sO(j,k))=0$
  whenever $-i\leq j+k\leq 0$, $-n\leq j\leq 0$ and $-m\leq k\leq 0$ but $(j,k)\not=(-n,0),(0,-m)$.
  \item[(b)] $E$ is a direct sum of the line bundles $\sO$, $\sO(0,1)$ and $\sO(1,0)$ with some balanced twist $(t,t)$.
  \end{enumerate}
  \end{theorem}
 Theorem \ref{t2} extends to the case of arbitrary $n, m$
  the classification of the ACM bundles on $\mathbb {P}^1\times \mathbb {P}^1$ proved in \cite{Kn}.
  In Theorem \ref{t0} we also prove an extension of Evans-Griffith criterion on $\mathbb {P}^n\times \mathbb {P}^m$.
  For a rank $r$ ($r<m+n$) vector bundle $E$ we demand the vanishing in ($a$) of the above theorems only for $i=1, \dots, r-1$ and we add a few extra cohomological vanishing condition in order to show that $E$ splits.  These extra  conditions do not appear in the Evans-Griffith criterion on $\mathbb P^{n}$.  We show that on $\mathbb {P}^n\times \mathbb {P}^m$  these extra hypotheses are necessary, and some conditions correspond to a direct summand $\sO\boxtimes\Omega_{\mathbf P^{m}}^a(a+1)$ (where $1\leq a\leq m-1$) or $\Omega_{\mathbf P^{n}}^a(a+1)\boxtimes\sO$ (where $1\leq a\leq n-1$). These bundles are also characterized  in \cite{cm0} by using a Beilinson's type spectral sequence.

  At the end of the paper we generalize the definition of regularity and the main results to  arbitrary multiprojective spaces.
  We work over an algebraically
closed field with characteristic
zero. The characteristic zero assumption is required for the proof
of \ref{t0},
 which use Le Potier vanishing theorem.

We thank E. Arrondo, A. P. Rao, and the referee for fundamental observations.

   \section{Regularity  on $\mathbf P^{n}\times\mathbf P^{m}$}\label{S2}

Let us consider $X=\mathbf P^{n}\times\mathbf P^{m}$. We want prove that if a  coherent sheaf $F$ is  regular according to Definition \ref{d1} then it is globally generated, and $F(p,p')$ is regular for $p,p'\geq 0$. We need the following Lemma:

\begin{lemma}Let $H$ be a generic hyperplane of $\mathbf P^{n}$. If $F$ is a regular coherent sheaf on $X$,
then $F_{|L_1}$ is  regular on $L_1=H\times \mathbf P^{m}$.\\
The similar statement is true for a generic hyperplane of $\mathbf P^{m}$.
\end{lemma}
\begin{proof}We follow the proof of \cite[(2.6.)]{hw}. We get this exact cohomology sequence:
$$\dots \rightarrow H^i(F(j,k)) \rightarrow H^i(F_{|L_1}(j,k)) \rightarrow H^{i+1}(F(j-1,k))\rightarrow \dots$$
If $j+k=-i$, $-n\leq j\leq 0$ and $-m\leq k\leq 0$, we also have $-n-1\leq j-1\leq 0$, so the first and the
third groups vanish  by hypothesis. Therefore the middle group vanishes as well. Hence $F_{|L_1}$ is regular.
\end{proof}

\begin{proposition}\label{p1} Let $F$ be a regular coherent sheaf on $X$.
  \begin{itemize}
  \item[(i)] $F(p,p')$ is regular for $p,p'\geq 0$.
  \item [(ii)] $H^0(F(k,k'))$ is spanned by $H^0(F(k-1,k'))\otimes H^0(\sO(1,0))$ if $k-1, k'\geq 0$, and by
  $H^0(F(k,k'-1))\otimes H^0(\sO(0,1))$ if $k, k'-1\geq 0$.
  \end{itemize}
  \end{proposition}

  \begin{proof} (i): We will prove part (i) by induction on $n+m$. Notice that for $n=0$, $X\cong \mathbf{P^m}$, and our definition of regularity coincides with
the Castelnuovo-Mumford regularity on $\mathbf{P^m}$. Hence we get part (i).  Let $F$ be a regular coherent sheaf on $X$. In
order to prove that  $F(1,0)$ is regular we follow the proof of \cite[(2.7.)]{hw}.\\
   Consider the exact
   cohomology sequence:
$$\dots \rightarrow H^i(F(j,k)) \rightarrow H^i(F(j+1,k)) \rightarrow H^{i}(F_{|L_1}(j+1,k))\rightarrow \dots$$
If $j+k=-i$, $-n\leq j\leq 0$ and $-m\leq k\leq 0$,  the first  group vanishes  by hypothesis. We want to show that the third group vanishes. Since $F_{|L_1}$ is regular by the above lemma and $F_{|L_1}(1,0)$ is regular by the inductive hypothesis, we have $H^i(F_{|L_1}(j+1,k))=0$
for $j+k=-i$, $-n+1\leq j\leq 0$ and $-m\leq k\leq 0$. Moreover the regularity of $F_{|L_1}(0,1)$ implies $H^i(F_{|L_1}(j+1,k))=0$ for $j=-n$, $k=n-i$ and $-m+1\leq k\leq 1$. If $(k,i,j) = (-m,n+m,-n)$, then we can use Grothendieck's vanishing, because
$dim (L_1) < n+m$. Therefore we have $H^i(F_{|L_1}(j,k))=0$
whenever $j+k=-i$, $-n\leq j\leq 0$ and $-m\leq k\leq 0$.
Therefore if $j+k=-i$, $-n\leq j\leq 0$ and $-m\leq k\leq 0$, then the middle group in the above sequence vanishes. This means that $F(1,0)$ is regular. A similar argument  works for $F(0,1)$.\\
(ii): We follow the  proof of \cite[(2.8.)]{hw}. We consider the following diagram:
$$  \begin{array}{ccc}
 H^0(F(k-1, k'))\otimes H^0(\sO(1,0))& \xrightarrow{\sigma}& H^0(F_{|L_1}(k-1,k'))\otimes
 H^0(\sO_{L_1}(1,0))\\
 \downarrow \mu& &\downarrow \tau\\
 H^0(F(k,k'))&\xrightarrow{\nu}& H^0(F_{|L_1}(k, k'))

 \end{array}$$
  Note that $\sigma$ is surjective if $k-1, k'\geq 0$  by the regularity condition. The map $\tau$ is surjective by (ii) applied to
$F_{|L_1}$.
  Since both $\sigma$ and $\tau$ are surjective, we can see, as in \cite{m} page $100$, that $\mu$ is
  surjective.
  \end{proof}

  \begin{remark}\label{gg} If  $F$ is a regular coherent sheaf on $X$, then it is globally generated.\\
  In fact by the above proposition we have the following surjections:
  $$H^0(F)\otimes
H^0(\sO(1,0))\otimes H^0(\sO(0,1))\rightarrow H^0(F(1,0))\otimes H^0(\sO(0,1)) \rightarrow H^0(F(1,1)),$$ so
the map
$H^0(F)\otimes H^0(\sO(1,1))\rightarrow H^0(F(1,1))$ is surjective.
  Moreover, we can consider a sufficiently large twist $l$ such that $F(l,l)$ is globally generated. The
  commutativity of the diagram
$$\begin{matrix}
H^0(F)\otimes H^0(\sO(l,l))\otimes\sO
&\to&H^0(F(l,l))\otimes\sO\\
\downarrow&&\downarrow\\
H^0(F)\otimes\sO(l,l)&\to&F(l,l)
\end{matrix}$$
yields the surjectivity of the map $H^0(F)\otimes\sO(l,l)\to F(l,l)$, which implies that $F$ is generated by its
global sections. \end{remark}
  \begin{remark}\label{b1} K\"{u}nneth's formula tells that $\sO(a,b)$ is regular if and only if $a \ge 0$ and $b \ge 0$.\\
  In fact $$H^{n+m}(\sO(a-n,b-m))\cong H^n(\sO(a-n))\otimes H^m(\sO(b-m))=0$$
  if and only if $a \ge 0$ or $b \ge 0$. Since $H^{n}(\sO(a-n,b))\cong H^n(\sO(a-n))\otimes H^0(\sO(b)),$ we see  that, if $\sO(a,b)$  is regular, we must have $a \ge 0$ and $b \ge 0$.\\
  In particular $\sO$ is regular while $\sO(-1,-1)$ is not and so $Reg(\sO)=0$.\\
  Moreover in a similar way we can see that $Reg(\sO\boxtimes\Omega_{\mathbf P^{m}}^a(a+1))=0$ for any $1\leq a\leq m-1$  \end{remark}
  Now we want to compare the two definitions of regularity.
  \begin{proposition}Let $F$ be a coherent sheaf on $X$.
  \begin{enumerate}
  \item[(i)] If $F$ is HW-regular, then it is regular.
  \item[(ii)] If $F$ is $(-m+1,-n+1)$-regular,  then it is HW-regular.
  \end{enumerate}
  \end{proposition}
  \begin{proof} (i): Let $F$ be HW-regular: for all $i>0$, $H^i(F\otimes
\sO(j,k))=0$
whenever $j+k=-i-1$, $j< 0$ and $ k< 0$.
By \cite[(2.7.)]{hw}, $F(p,p')$ is HW-regular for any $p\geq 0$ and $p'\geq 0$.
In particular $H^i(F\otimes
\sO(j,k))=0$
whenever $i>0$, $j+k=-i$, $-n\leq j\leq 0$ and $-m\leq k\leq 0$.

(ii): Let $F$ be $(-m+1,-n+1)$-regular. We have
$H^{m+n}(F(-n-m+1,-m-n+1))=0,$ and so
$H^{m+n}(F(-n-m+1,-1))=H^{m+n}(F(-n-m+2,-2))=\dots =H^{m+n}(F(-1,-m-n+1))=0.$
In the same way all  other vanishing conditions in the definition of regularity by Hoffmann and Wang are satisfied.
\end{proof}

  \section{Splitting Criteria for Vector Bundles}\label{S3}
  We use our notion of regularity in order to  prove some splitting criteria on $X=\mathbf P^{n}\times\mathbf P^{m}$.\\

  We need the following definition:
  \begin{definition} We say that a vector bundle $E$ on $X$ is ACM if, for any $i=1, \dots, m+n-1$ and for any integer $t$,  $$H^i(E(t,t))=0.$$
\end{definition}
\begin{remark}\label{b2}K\"{u}nneth's formula implies that $\sO(a,b)$ is ACM if and only if $a-b \ge -m$ and $b-a \ge -n$.\\
We recall that on $X$ we have the following Koszul complexes:\\
 \smallbreak \centerline{${}\hfill 0\to \sO(-n-1,-m-1)\to \sO(-n,-m-1)^{\binom{n+1}{n}}\to \dots \to \sO(0,-m-1)\to 0,$ \hfill $[K_1]$}
 \smallbreak  \centerline{${}\hfill 0\to \sO(0,-m-1)\to \sO(0,-m)^{\binom{m+1}{m}}\to \dots \to \sO\to 0,$ \hfill $[K_2]$}
 \smallbreak
  \centerline{${}\hfill 0\to \sO(-n-1,-m-1)\to \sO(-n,-m-1)^{\binom{n+1}{n}}\to\dots \to \sO(0,-m-1)\to \dots \to \sO\to 0$ \hfill $[K_3]$}\smallbreak
From the above exact sequences in cohomology we get the following isomorphisms:
$$\textrm{$\epsilon_1: H^m(\sO(0,-m-1))\to H^{n+m}(\sO(-n-1,-m-1))$ and $\epsilon_2: H^0(\sO)\to H^{m}(\sO(0,-m-1))$}$$
  \end{remark}
  Now are ready to prove Theorem \ref{t1}:\\

  {\emph {Proof of Theorem \ref{t1}.} $(i) \Rightarrow (ii)$. Let us assume that $t$ is an integer such that $E(t,t)$ is regular, but
  $E(t-1,t-1)$ is not regular.\\
By the definition of regularity and (i) we can say that $E(t-1,t-1)$ is not regular if and only if
$H^{m+n}(E(t-1,t-1)\otimes \sO(-n,-m))\not=0$. Serre duality gives $H^0(E^{\vee}(-t,-t))\not=0$.
Since $E(t,t)$ is globally generated by Remark \ref{gg}, and $H^0(E^{\vee}(-t,-t))\not=0$, we conclude
 that $\sO$ is a direct summand of $E(t,t)$.   By iterating these arguments we get (ii).\\
      $(ii)\Rightarrow (i)$. $\sO(j,k)$ is ACM for $-n\leq j\leq 0$ and $-m\leq k\leq 0$. Hence $E\cong \bigoplus_{i=1}^r \sO(t_i,t_i)$ satisfies all the conditions in (i).\\
      \qquad

Now we  prove Theorem \ref{t2}:\\

  {\emph {Proof of Theorem \ref{t2}.}} $(a)\Rightarrow (b)$. By Serre duality
  $H^i(E(t,t)\otimes \sO(j,k))\cong H^{m+n-i}(E^{\vee}(t,t)\otimes \sO(-n-1-j,-m-1-k)).$
  Let us assume that $t$ is an integer such that $E(t,t)$ is regular, but
  $E(t-1,t-1)$ is not regular.
By the definition of regularity and (a) we can say that $E(t-1,t-1)$ is not regular if and only if one of the
following conditions holds:
     \begin{enumerate}
     \item[(i)] $H^{m+n}(E(t-1,t-1)\otimes \sO(-n,-m))\not=0$,
      \item[(ii)] $H^n(E(t-1,t-1)\otimes \sO(-n,0))\not=0$.
   \item[(iii)] $H^m(F(t-1,t-1)\otimes \sO(0,-m))\not=0$. \end{enumerate}
  Let us consider the obove conditions one by one:\\
     (i): Let $H^{m+n}(E(t-1,t-1)\otimes \sO(-n,-m))\not=0$. Then we conclude that $\sO(t,t)$ is a direct
summand as in the above proof.\\
(ii): Let $H^n(E(t-1,t-1)\otimes \sO(-n,0))\not=0$. Let us consider the exact sequence $[K_1]$ tensored with $E(t,t+m)$:

      $$0 \to \sO(-n-1,-1)\otimes E(t,t) \rightarrow \sO(-n,-1)^{\binom{n+1}{n}}\otimes E(t,t)\rightarrow \dots$$
       $$\dots\rightarrow\sO(-1,-1)^{\binom{n+1}{1}}\otimes E(t,t) \rightarrow \sO(0,-1)\otimes E(t,t)
\to 0.$$
Since $$ H^n(E(t-n,t-1))= \dots =H^1(E(t-1,t-1))=0,$$
(if $n>m$ we also  use the hypothesis $H^m(E(t-m,t-1))=0$)  we have a surjective map
 $$\alpha_1: H^0(E(t,t-1))\to H^{n}(E(t-n-1,t-1)).$$
On the other hand  $$H^{n}(E(t-n-1,t-1))\cong H^m(E^{\vee}(-t,-t-m))$$ so let us consider sequence $[K_2]$ tensored with $E^\vee(t,t+1)$
 $$0 \to \sO(0,-m)\otimes E^{\vee}(-t,-t) \rightarrow \sO(0,-m+1)^{\binom{m+1}{m}}\otimes
  E^{\vee}(-t,-t)\rightarrow \dots $$ $$\dots\rightarrow\sO(0,0)^{\binom{m+1}{1}}\otimes
  E^{\vee}(-t,-t) \rightarrow \sO(0,1)\otimes E^{\vee}(-t,-t)
\to 0.$$
Since
$$H^n(E(t-n-1,t-2))=\dots =H^{m+n-1}( E(t-n-1,t-m-1))=0,$$
(if $m>n$ we also use the hypothesis $H^m(E(t-n-1,t-2-m+n))=0$) and by Serre duality $$H^m(E^{\vee}(-t,-t-m+1))=\dots =H^1( E^{\vee}(-t,-t))=0,$$
 we have a surjective map  $$\alpha_2: H^0(E^{\vee}(-t,-t+1))\to H^{m}(Eì\vee(-t,-t-m)).$$
 Let us consider the following diagram:
 $$\begin{array}{ccc} H^{n}(E(t-n-1,t-1))\otimes H^m(E^{\vee}(-t,-t-m))& \xrightarrow{\sigma}& H^{m+n}(\sO(-n-1,-1)\otimes \sO(0,-m))\\ \uparrow\alpha_1\otimes 1 & &\uparrow\epsilon'_1 \\
 H^0(E(t,t-1)) \otimes H^m(E^{\vee}(-t,-t-m))& \xrightarrow{\mu}& H^{m}(\sO(0,-1)\otimes \sO(0,-m))\\
 \uparrow1\otimes\alpha_2 & &\uparrow\epsilon'_2 \\
      H^0(E(t,t-1))  \otimes H^0(E^{\vee}(-t,-t+1))& \xrightarrow{\tau}& H^{0}(\sO(0,-1)\otimes \sO(0,1)),\\
        \end{array}$$
        where $\epsilon'_1$ and $\epsilon'_2$ are naturally induced by $\epsilon_1$ and $\epsilon_2$ (see remark \ref{b2}) and the horizontal maps $\sigma, \mu , \tau$
        are the Yoneda pairing surjections inducing Serre duality, i.e., $ \sigma$ induces the duality between
$H^n(E(t-n-1,t-1))$ and $H^m(E^\vee (-t,-t-m))$ (and similarly for $\mu$ and $\tau$) (see \cite{ak}, p. 67 and Theorem at p. 1,
\cite{h}, Theorem III.7.6 (a)).
        The maps $\alpha _i$ and $\epsilon _i$ are obtained from the same Koszul complex  by splitting it into short exact sequences
        and then applying some vanishing. Yoneda pairing commutes with the connecting morphisms induced by short exact sequences (see \cite{ak}, Theorem (1.1) at p. 67).
        Hence  the horizontal maps $\sigma, \mu, \tau$ commute with the connecting maps (the vertical maps $\alpha_1\otimes 1, 1\otimes\alpha_2, \epsilon'_1, \epsilon'_2$). Therefore the above diagram is commutative.\\
        Since $\sigma$ induces the Serre duality
isomorphism $H^{n}(E(t-n-1,t-1))\cong H^m(E^{\vee}(-t,-t-m))$, there are non-zero elements $s\in H^{n}(E(t-n-1,t-1))$ and $s'\in H^m(E^{\vee}(-t,-t-m))$ such that $\sigma(s\otimes s')\ne 0$. Since $\alpha_1$ and $\alpha_2$ are surjective, there are
 $g\in H^{0}(E(t,t-1))$ and $f\in H^0(E^{\vee}(-t,-t+1))$ such that $\tau(f,g) \ne 0$. The elements $f$ and $g$ can be regarded as elements of ${Hom}(E(t,t), \sO(0,1))$ and ${Hom}(\sO(0,1), E(t,t))$ respectively.    This means that the  map $$f\circ g : \sO(0,1) \to \sO(0,1)$$ is non-zero and hence it is an
isomorphism.
This isomorphism shows that $\sO(0,1)$ is a direct summand of $E(t,t)$.\\
(iii): Assume $H^m(F(t-1,t-1)\otimes \sO(0,-m))\ne 0$. Arguing as above we conclude
 that $\sO(1,0)$ is a direct summand of $E(t,t)$.\\
  $(b)\Rightarrow (a)$. As in Theorem \ref{t1}.\qed
\begin{remark}The case $n=m=1$ recovers the classification of  ACM bundles on $\mathbb {P}^1\times \mathbb {P}^1$ (see \cite{Kn}).
The proof in this case coincides with the one of \cite{bm2}, Theorem 1.2.
\end{remark}
 \begin{corollary}\label{r2}Let $E$ be a  vector bundle on $X$, $a$ and $b$ be two integers.\\ Then the following conditions are equivalent:
  \begin{itemize}
  \item[(i)] for any $i=1, \dots, m+n-1$ and any integer $t$,  $$H^i(E(a+t,b+t)\otimes \sO(j,k))=0$$
  whenever $-i\leq j+k\leq 0$, $-n\leq j\leq 0$ and $-m\leq k\leq 0$ but $(j,k)\not=(-n,0),(0,-m)$.\\

  \item[(ii)] $E$ is a direct sum of  the line bundles $\sO(a,b)$, $\sO(a,b+1)$ and $\sO(a+1,b)$ with some balanced twist $(t,t)$.
  \end{itemize}
  \end{corollary}\begin{proof} Consider $E\otimes \sO(-a,-b)$ and  apply the previous theorem.
  \end{proof}
By applying Le Potier  vanishing Theorem (see \cite{le} or \cite[(7.3.5.)]{LA}) we can prove the following result:
\begin{theorem}\label{t0}Let $E$ be a rank $r$  vector bundle on $\mathbb {P}^n\times \mathbb {P}^m$ with $Reg (E)= 0$.\\ Then the following conditions are equivalent:
  \begin{enumerate}
  \item[(1)] for any $i=1, \dots, min(r,m+n)-1$,  $$H^i(E(-1,-1)\otimes \sO(j,k))=0$$ whenever
  $j+k\geq -i$, $-n< j\leq 0$ and $-m< k\leq 0$.\\
  Moreover for any $u=1,\dots ,m-1$, $H^{m+u}(E(-1,-1)\otimes \sO(-n,-u-1))=0$ and for any $v=1,\dots ,n-1$, $H^{n+v}(E(-1,-1)\otimes \sO(-v-1,-m))=0$.
\item[(2)] $E$ has one of the following bundles as a direct summand: $\sO$, $\sO(0,1)$, $\sO(1,0)$, $\sO\boxtimes\Omega_{\mathbf P^{m}}^a(a+1)$ (where $1\leq a\leq m-1$) or $\Omega_{\mathbf P^{n}}^a(a+1)\boxtimes\sO$ (where $1\leq a\leq n-1$).
  \end{enumerate}
  \end{theorem}
 \begin{proof} $(1)\Rightarrow (2)$. Since $Reg (E)= 0$, the bundle $E$ is regular, while $E(-1,-1)$ is not.\\
  $E$ is globally generated by Remark \ref{gg}. Since the tensor product of a spanned vector bundle by an ample vector
bundle is ample (see \cite{ha} Corollary III.1.9), we have
  $$\textrm{$a, b>0$ $\Rightarrow E(a,b)$ is ample. }$$
  Let us assume $r<m+n$. Therefore Le Potier vanishing theorem gives $H^i(E^{\vee}(-a,-b))=0$ for every $a,b>0$ and $i=1, \dots, n+m-r$.\\
 So by Serre duality $H^i(E(-n-1+a,-m-1+b))=0$ for every $a,b>0$ and $i=r, \dots, n+m-1$.\\
 In particular for any $i=1, \dots, m+n-1$,  $H^i(E(-1,-1)\otimes \sO(j,k))=0$ whenever
  $j+k\geq -i$, $-n< j\leq 0$ and $-m< k\leq 0$. We can then say that $E(-1,-1)$ is not regular if and only if one of the
following conditions is satisfied (if $r\geq m+n$ we do not need Le Potier vanishing theorem):
     \begin{enumerate}

   \item [(i)] there exists an integer $a\in \{0,\dots ,m\}$ such that $H^{n+a}(E(-1,-1)\otimes \sO(-n,-a))\ne 0$,
   \item [(ii)] there exists an integer $a\in \{0,\dots ,n\}$ such that $H^{m+a}(E(-1,-1)\otimes \sO(-a,-m))\ne 0$. \end{enumerate}
(i) If $a=m$ or $a=0$, by the proof of Theorem \ref{t2} we have direct summands $\sO$ or $\sO(0,1)$ respectively.\\
 Fix $a\in \{1,\dots m-1\}$ and assume $H^{n+a}(E(-1,-1)\otimes \sO(-n,-a))\ne 0$. Let us consider the exact sequence $[K_1]$ twisted by $(0,-1-a)$
 and the dual of
\smallbreak
  \centerline{${}\hfill 0 \to \sO\boxtimes \Omega_{\mathbf P^{m}}^a(a+1) \rightarrow \sO(0,1)^{\binom{m+1}{a}}\rightarrow \sO(0,2)^{\binom{m+1}{a-1}}\to\dots$\hfill}
  \centerline{${}\hfill\dots\rightarrow\sO(0,a)^{\binom{m+1}{1}} \rightarrow \sO(0,a+1)
\to 0.$\hfill $[K_4]$}
We tensor by $E$ and obtain
      $$0 \to \sO(-n-1,-1-a)\otimes E \rightarrow \sO(-n,-1-a)^{\binom{n+1}{n}}\otimes E\rightarrow \dots$$
       $$\dots\rightarrow\sO(-1,-1-a)^{\binom{n+1}{1}}\otimes E \rightarrow \sO(0,-a)^{\binom{m+1}{1}}\otimes E\rightarrow\dots\rightarrow \sO\boxtimes(\Omega_{\mathbf P^{m}}^a)^\vee(-1-a)\otimes E
\to 0,$$
Since $$ H^{n+a}(E(-n,-1-a))= \dots =H^{a+1}(E(-1,-1-a)) =H^{a}(E(0,-a))=\dots =H^1(E(0,-1))=0,$$
  we have a surjective map
 $$\alpha_1: H^0(\sO\boxtimes(\Omega_{\mathbf P^{m}}^a)^\vee(-1-a)\otimes E)\to H^{n+a}(E(-n-1,-1-a)).$$
On the other hand  $$H^{n+a}(E(-n-1,-a-1))\cong H^{m-a}(E^{\vee}(0,-m+a)),$$ so let us consider the exact sequence
 $$\sO(0,-m+a)\rightarrow\sO(0,-m+a+1)^{\binom{m+1}{m}}\rightarrow\dots\rightarrow\sO^{\binom{m+1}{a+1}} \rightarrow \sO\boxtimes\Omega^a_{\mathbf P^{m}}(a+1)
\to 0.$$
Tensoring by $E^\vee$, $$0 \to E^\vee(0,-m+a) \rightarrow E^\vee(0,-m+a+1)^{\binom{m+1}{m}}\rightarrow E^\vee(0,-m+a+2)^{\binom{m+1}{m-1}}\rightarrow \dots $$ $$\dots\rightarrow E^\vee(0,0)^{\binom{m+1}{a+1}} \rightarrow \sO\boxtimes\Omega^a_{\mathbf P^{m}}(a+1)\otimes E^\vee
\to 0.$$
 Since $$H^{m+n-1}(E(-n-1,-m-1))=\dots =H^{m+n-m+a+1}( E(-n-1,-a-3))=$$ $$=H^{m+n-m+a}(E(-n-1,-a-2))=0,$$ and by Serre duality
 $$H^{1}(E^{\vee})=\dots =H^{m-a-1}( E^{\vee}(0,-m+a+2))=H^{m-a}(E^{\vee}(0,-m+a+1))=0,$$
 we have a surjective map  $$\alpha_2: H^0(\sO\boxtimes\Omega^a_{\mathbf P^{m}}(a+1)\otimes E^{\vee})\to H^{m-a}(E^{\vee}(0,-m+a)).$$
 From the surjective maps  $\alpha _i$, $i=1,2$, we get surjections
$$\alpha _1\otimes 1: H^0(\sO\boxtimes(\Omega_{\mathbf P^{m}}^a)^\vee (-1-a)\otimes E)\otimes
H^{m-a}(E^{\vee}(0,-m+a))\to$$ $$\to H^{n+a}(E(-n-1,-1-a))\otimes
H^{m-a}(E^{\vee}(0,-m+a)),$$
$$1\otimes \alpha _2: H^{n+a}(E(-n-1,-1-a))\otimes H^0(\sO\boxtimes\Omega^a_{\mathbf P^{m}}(a+1)\otimes E^{\vee}) \to$$ $$\to H^{n+a}(E(-n-1,-1-a))\otimes H^{m-a}(E^{\vee}(0,-m+a)).$$
We use the Yoneda pairings
$$\tau : H^0(\sO\boxtimes(\Omega_{\mathbf P^{m}}^a)^\vee(-1-a)\otimes E)\otimes H^0(\sO\boxtimes\Omega^a_{\mathbf P^{m}}(a+1)\otimes  E^{\vee}) \to H^0(\sO ) \cong \mathbb {C},$$
$$\mu :H^0(\sO\boxtimes(\Omega_{\mathbf P^{m}}^a)^\vee(-1-a)\otimes E)\otimes H^{m-a}(E^{\vee}(0,-m+a)) \to H^{m-a}(\sO\boxtimes\Omega^a_{\mathbf P^{m}}(a+1)) \cong \mathbb {C},$$
\begin{align*}
&\sigma : H^{n+a}(\sO\boxtimes(\Omega_{\mathbf P^{m}}^a)^\vee(-1-a)\otimes E(-n-1,-m))\otimes
H^{m-a}(\sO\boxtimes\Omega^a_{\mathbf P^{m}}(a+1)\otimes E^{\vee}(0,-m+a)) \\
&\to H^{n+m}(\sO (-n-1,-m-1) \cong \mathbb {C}.\end{align*}The exact sequence
$[K_4]$ gives surjections $\epsilon _1:  H^{m-a}(\sO\boxtimes(\Omega_{\mathbf P^{m}}^a)^\vee (-1-a))\to H^{n+m}(
\sO(-n-1,-m-1)) \cong \mathbb {C}$ and $\epsilon _2:  H^0(\sO)
\to H^{m-a}(\sO\boxtimes\Omega^a_{\mathbf P^{m}}(a+1))) \cong \mathbb {C}$.
 As in the proof of Theorem \ref{t2}, from the surjective maps $\alpha_1, \alpha_2$  and the non-degenerate pairings $\sigma, \tau$ we can find
  $$\textrm{ $ g: \sO\boxtimes\Omega^a_{\mathbf P^{m}}(a+1)\rightarrow  E$ and $f: E \rightarrow  \sO\boxtimes\Omega^a_{\mathbf P^{m}}(a+1)$}$$ such that $f\circ g$ is a nonzero map. Since $\Omega^a_{\mathbf P^{m}}(a+1)$ is
simple,  we conclude that $\sO\boxtimes\Omega^a_{\mathbf P^{m}}(a+1)$ is a direct summand of $E$.\\
 (ii): Arguing as above we get that $E$ has one of the following bundles as a direct summand: $\sO$, $\sO(1,0)$,  or $\Omega_{\mathbf P^{n}}^a(a+1)\boxtimes\sO$ (where $1\leq a\leq n-1$).\\
  $(2)\Rightarrow (1)$. We have to check that for any $a=1,\dots, m-1$, $\sO\boxtimes\Omega^a_{\mathbf P^{m}}(a+1)$ satisfies all the conditions of $(1)$.
  Let us consider all the groups of cohomology that can be different from zero:\\
  $H^a(\sO(j)\boxtimes\Omega^a_{\mathbf P^{m}}(a+1+k))\not=0$ if and only if $j\geq 0$ and $k=-a-1$,\\
  $H^n(\sO(j)\boxtimes\Omega^a_{\mathbf P^{m}}(a+1+k))\not=0$ if and only if $j\leq -n-1$ and $k\geq -1$,\\
  $H^m(\sO(j)\boxtimes\Omega^a_{\mathbf P^{m}}(a+1+k))\not=0$ if and only if $j\geq 0$ and $k\leq -m-a-1$, and\\
  $H^{n+a}(\sO(j)\boxtimes\Omega^a_{\mathbf P^{m}}(a+1+k))\not=0$ if and only if $j\leq -n-1$ and $k=-a-1$.\\
  So the conditions $(1)$ are satisfied.
 \end{proof}

We can easily generalize the notion of regularity on $X=\mathbf P^{n_1}\times\dots \times\mathbf P^{n_s}$
      ($d=n_1+\dots+n_s$) and adapt Theorems \ref{t1}, and \ref{t2},  to this setup. The proofs  are very similar to the  case $s=2$. Therefore we do not include them here.
      \begin{definition}\label{d4}
A  coherent sheaf $F$ on $X=\mathbf P^{n_1}\times\dots \times\mathbf P^{n_s}$  is said to be {\it $(p_1, \dots,
p_s)$-regular} if, for all $i>0$,
$$H^i(F(p_1, \dots, p_s)\otimes \sO(k_1, \dots, k_s))=0$$ whenever $k_1+ \dots, +k_s=-i$ and $-n_j\leq k_j\leq 0$ for any
$j=1,\dots , s$.\\
\end{definition}
  \begin{remark}\label{b3}K\"{u}nneth's formula implies that $\sO(a_1,\dots, a_s)$ is ACM if and only if for any $j=1,\dots , s$ there are $h,k\not=j$ such that $a_j-a_h\leq n_h$ and $a_j-a_k\geq -n_j$.\\
  \end{remark}
  \begin{proposition}\label{p1b} Let $F$ be a regular coherent sheaf on $X$. Then
  \begin{enumerate}
  \item $F(p_1, \dots ,p_s)$ is regular for $p_1, \dots ,p_s\geq 0$.
  \item For any $j=1,\dots, s$, $H^0(F(k_1, \dots, k_s))$ is spanned by
   $$H^0(F(k_1, \dots, k_j-1, \dots, k_s))\otimes H^0(\sO(0,\dots ,1,\dots, 0))$$ if
   $k_1, \dots, k_j-1, \dots, k_s\geq 0$.
   \item $F$ is globally generated
  \end{enumerate}
  \end{proposition}
  We can now give the following splitting criteria which are the generalizations of Theorem \ref{t1}, \ref{t2}:
\begin{theorem}\label{t3}Let $E$ be a rank $r$ vector bundle on $X=\mathbf P^{n_1}\times\dots \times\mathbf P^{n_s}$.
Set $d=n_1+\dots+n_s$. Then the following conditions are equivalent:
  \begin{enumerate}
  \item for any $i=1, \dots, d-1$ and for any integer $t$,  $H^i(E(t,\dots, t)\otimes \sO(k_1, \dots, k_s))=0$
   whenever $k_1+ \dots, +k_s=-i$ and $-n_j\leq k_j\leq 0$ for any $j=1,\dots , s$.
  \item There are $r$ integer $t_1, \dots, t_r$ such that $E\cong \bigoplus_{i=1}^r \sO(t_i,\dots, t_i)$.
  \end{enumerate}
  \end{theorem}
\begin{proof}  $(1)\Rightarrow (2)$. Let assume that $t$ is an integer such that $E(t, \dots, t)$ is regular but
  $E(t-1,\dots, t-1)$ not.
By the definition of regularity and $(1)$ we can say that $E(t-1,\dots, t-1)$ is not regular if and only if
$H^{d}(E(t-1, \dots, t-1)\otimes \sO(-n_1,\dots, -n_s))\not=0$.
 Now since $E(t,\dots, t)$ is globally generated  and $H^0(E^{\vee}(-t,\dots, -t))\not=0$ we can conclude
 that $\sO$ is a direct summand of $E(t,\dots, t)$.
      By iterating these arguments we get $(2)$.\\
      $(2)\Rightarrow (1)$. $E\cong \bigoplus_{i=1}^r \sO(t_i, \dots, t_i)$ is ACM then it satisfies all the conditions in $(1)$.\end{proof}
      \begin{theorem}\label{t2b}Let $E$ be a rank $r$ vector bundle on $X=\mathbf P^{n_1}\times\dots \times\mathbf P^{n_s}$.
Set $d=n_1+\dots+n_s$. Then the following conditions are equivalent:
  \begin{enumerate}
  \item for any $i=1, \dots, d-1$ and for any integer $t$,  $H^i(E(t,\dots, t)\otimes \sO(k_1, \dots, k_s))=0$
   whenever $k_1+ \dots, +k_s\geq -i$ and $-n_j-1\leq k_j\leq 1$ for any $j=1,\dots , s$,  but with an index $j$ such that $k_j\not=0, -n_j$.
  \item $E$ is a direct sum of  line bundles  $\sO(l_1,\dots , l_s)$ (where  $l_j=1$ or $l_j=0$ for any $j=1,\dots , s$)  with some balanced twist $(t,\dots , t)$.
  \end{enumerate}
  \end{theorem}
  \begin{proof}
  $(1)\Rightarrow (2)$.
  Let assume that $t$ is an integer such that $E(t,\dots, t)$ is regular but
  $E(t-1,\dots, t-1)$ not.
By the definition of regularity and $(1)$ we can say that $E(t-1,\dots, t-1)$ is not regular if and only if one of the
following conditions is satisfied:
     \begin{enumerate}
     \item[i] $H^{d}(E(t-1, \dots, t-1)\otimes \sO(-n_1,\dots, -n_s))\not=0$,
      \item[ii] there are $s$ numbers $h_1, \dots, h_s$ where for any $j=1,\dots , s$ $h_j=0$ or $h_j=n_j$ and  $0< k_1, \dots, k_s< d$ such that\\ $H^{h_1+ \dots + h_s}(E(t-1,\dots, t-1)\otimes \sO(-h_1, \dots, -h_s))\not=0$.
    \end{enumerate}
  Let us consider one by one the conditions:\\
     $(i)$ Let $H^{d}(E(t-1, \dots, t-1)\otimes \sO(-n_1,\dots, -n_s))\not=0$, we can conclude that $\sO(t,\dots ,t)$ is a direct
summand as in the above theorem.\\
$(ii)$ We can assume that there is an integer $l$ with $1\leq l<s$ such that\\  $H^{n_1+\dots +n_l}(E(t-1,\dots ,t-1)\otimes \sO(-n_1,\dots ,-n_l,0, \dots, 0))\not=0$.\\
 Let us consider the following exact sequences tensored by $E(t, \dots, t)$:
 $$0 \to \sO(-n_1-1,\dots ,-n_l-1,-1, \dots, -1) \rightarrow \dots\rightarrow \sO(0,-n_2-1\dots ,-n_l-1,-1, \dots, -1)
\to 0,$$
$$0 \to \sO(0,-n_2-1,\dots ,-n_l-1,-1, \dots, -1) \rightarrow \dots\rightarrow \sO(0,0,-n_3-1\dots ,-n_l-1,-1, \dots, -1)
\to 0,$$
$$\dots$$
$$0 \to \sO(0,\dots ,0,-n_l-1,-1, \dots, -1) \rightarrow \dots\rightarrow \sO(0,\dots ,0,-1, \dots, -1)
\to 0.$$
Since, by  the vanishing conditions in $(1)$ $$H^{n_1+\dots +n_l}(E(t-n_1, t-n_2-1,\dots ,t-n_l-1, t-1,\dots ,t-1))=\dots$$ $$\dots =H^{1+n_2+\dots +n_l}(E(t-1,t-n_2-1,\dots ,t-n_l-1, t-1,\dots ,t-1))=$$ $$=H^{2+n_2+\dots +n_l}(E(t,t-n_2,\dots ,t-n_l-1, t-1,\dots ,t-1))=\dots $$ $$\dots=H^{1}(E(t,\dots ,t,t-1, t-1,\dots ,t-1))=0,$$

we can conclude that $$H^0(E(t,\dots, t)\otimes\sO(0,\dots , 0,-1,\dots ,-1))\not=0.$$
On the other hand  $$H^{n_1+\dots +n_l}(E(t-1,\dots ,t-1)\otimes \sO(-n_1,\dots ,-n_l,0, \dots, 0))\cong$$ $$\cong H^{n_{l+1}+\dots +n_s}(E^{\vee}(-t,\dots,-t)\otimes(0,\dots ,0,-n_{l+1}, \dots,-n_s)).$$

Let us consider the following exact sequences tensored by $E^{\vee}(-t, \dots, -t)$:

 $$0 \to \sO(0,\dots ,0,-n_{l+1}, \dots, -n_s) \rightarrow \dots\rightarrow \sO(0,\dots ,0,1,-n_{l+2} \dots, -n_s)
\to 0,$$
$$0 \to \sO(0,\dots ,0,1,-n_{l+2}, \dots, -n_s)) \rightarrow \dots\rightarrow \sO(0,\dots ,0,1,1,-n_{l+3} \dots, -n_s)
\to 0,$$
$$\dots$$
$$0 \to \sO(0,\dots ,0,1, \dots,1,-n_s) \rightarrow \dots\rightarrow \sO(0,\dots ,0,1, \dots, 1)
\to 0.$$
Since, by the vanishing conditions in $(1)$,
 $$H^{n_{1}+\dots +n_l}(E(t-n_1-1,\dots ,t-n_l-1,t-2, t-1, \dots, t-1))=\dots$$
 $$\dots=H^{n_1+\dots +n_l+n_{l+1}-1}(E(t-n_1-1,\dots ,t-n_l-1, t-n_{l+1}-1, t-1, \dots, t-1))=$$
 $$=H^{n_1+\dots +n_l+n_{l+1}-1}(E(t-n_1-1,\dots ,t-n_l-1, t-n_{l+1}-2, t-2, \dots, t-1))    =\dots $$
 $$\dots =H^{d-1}(E(t-n_1-1,\dots ,t-n_l-1,t-n_{l+1}-2,\dots , t-n_{s_1}-2, t-n_s-1))=0,$$

 and by Serre duality,

 $$H^{n_{l+1}+\dots +n_s}(E^\vee(-t,\dots ,-t,-t-n_{l+1}+1, -t-n_{l+2}, \dots, -t-n_s))=\dots$$
 $$\dots=H^{1+n_{l+2}+\dots +n_s}(E^\vee(-t,\dots ,-t,-t, -t-n_{l+2}, \dots, -t-n_s))=$$
 $$=H^{n_{l+2}+\dots +n_s}(E^\vee(-t,\dots ,-t,-t+1, -t-n_{l+2}+1, \dots, -t-n_s))=\dots $$
 $$\dots =H^{1}(E^\vee(-t,\dots ,-t,-t+1,\dots , -t+1, -t))=0$$

 we can conclude that $$H^0(E^\vee(-t,\dots, -t)\otimes\sO(0,\dots , 0,1,\dots ,1))\not=0.$$
So by arguing as in Theorem \ref{t2} we have that $\sO(0,\dots ,0,1, \dots, 1)$ is a direct summand of $E(t,\dots, t)$.
  \end{proof}
  \begin{remark}We could also generalize Theorem
\ref{t0}. We preferred
not to do it explicitly, since it represents a small improvement
compared with the difficulty to write it in a clear way.
\end{remark}

\bibliographystyle{amsplain}

\end{document}